\documentclass[12pt]{article}
\usepackage{amssymb,amsfonts,amsthm,amsmath,graphicx}

\newtheorem*{theorem}{Theorem}

\newtheorem*{definition}{Definition}

\newtheorem*{proposition}{Proposition}
\newtheorem{lemma}{Lemma}

\newcommand{\bl}{\begin{lemma}}
\newcommand{\el}{\end{lemma}}

\newcommand{\bp}{\begin{proof}}
\newcommand{\ep}{\end{proof}}

\newcommand{\R}{\mathbb R}
\newcommand{\eps}{\varepsilon}
\newcommand{\wt}{\widetilde}
\newcommand{\wh}{\widehat}

\newcommand{\mfrac}[2]{\scalebox{1.3}{\ensuremath{\frac{#1}{#2}}}}

\begin{document}
\title{The Newman algorithm for constructing polynomials with restricted coefficients and many real roots}
\author{Markus Jacob and Fedor Nazarov}
\maketitle

\begin{abstract}
Under certain natural sufficient conditions on the sequence of uniformly bounded closed sets $E_k\subset\R$ of admissible coefficients, we construct a polynomial $P_n(x)=1+\sum_{k=1}^n\eps_k x^k$, $\eps_k\in E_k$, with at least $c\sqrt n$ distinct roots in $[0,1]$, which matches the classical upper bound up to the value of the constant $c>0$. Our sufficient conditions cover the Littlewood ($E_k=\{-1,1\}$) and Newman ($E_k=\{0,(-1)^k\}$) polynomials and are also necessary for the existence of such polynomials with arbitrarily many roots in the case when the sequence $E_k$ is periodic.
\end{abstract}

\section{Acknowledgements}
\label{acknowledgements}

This article arose from the MathOverflow discussion \cite{MOP}. We thank all the discussion participants for their contributions and MathOverflow itself for providing the online space and friendly environment for the exchange of ideas. In addition, our special thanks go to Tamas Erd\'elyi, who attracted our attention to the Newman Decomposition Lemma, and to Mikhail Sodin, who read the first draft of this article and made numerous suggestions that greatly improved its readability. The second named author was partially supported by the NSF grant DMS-1900008.

\section{Introduction}
\label{intro}

Consider a sequence $\mathcal E$ of uniformly bounded closed sets $E_1,E_2,\dots,E_k,\ldots\subset\R$. The uniform boundedness condition means that there exists $A\in(0,+\infty)$ such that $E_k\subset[-A,A]$ for all $k\in\mathbb N$. Let $\mathcal P=\mathcal P(\mathcal E)$ be the set of all polynomials of the form
$$
P_n(x)=1+\sum_{k=1}^n\eps_k x^k,\qquad \eps_k\in E_k\,.
$$
What is the maximal number $r$ of distinct roots a polynomial $P_n\in\mathcal P$ of a given degree $n$ can have in $[0,1]$ or, equivalently, for a given $r\in\mathbb N$, what is the lowest possible degree $n$ of a polynomial $P_n\in\mathcal P$ with at least $r$ distinct roots in $[0,1]$?

The now classical result is that the inequality 
\begin{equation}
r\le C(A)\sqrt n
\label{r2len}
\end{equation}
always holds with some $C(A)\in(0,+\infty)$ depending on $A$ only. One can find this estimate, for instance, in \cite{BEK}, which goes back to 1999. For reader's convenience, we also present a proof in the appendix.

The question then becomes if (or when) this bound is asymptotically sharp for large $n$ up to the value of the numerical constant $C(A)$. As far as we know, in this form it has not been previously answered even for such natural and well-known families as Littlewood polynomials ($E_k=\{-1,1\}$) and Newman polynomials\footnote{The standard definition of the Newman polynomials restricts the coefficients to the set $\{0,1\}$ but then the question about roots should be asked on $[-1,0)$, so we took the liberty to change the variable to $-x$ to place this family into our general framework.}
($E_k=\{0,(-1)^k\}$). In both cases the best published lower bounds for $r$ seem to be polylogarithmic in terms of $n$.  Some more substantial progress has been made for the polynomials with integer coefficients of height $1$ ($E_k=\{-1,0,1\}$). In this case it was shown in \cite{E} that $r\ge cn^{1/4}$. The interested reader can find more related results and the general overview of the history of the question in \cite{BEK}. 

To motivate our next definition, let us consider one simple obstacle that prevents polynomials in $\mathcal P(\mathcal E)$ from having many roots in $[0,1]$ regardless of their degree. Suppose that the sequence $E_k$ is $M$-periodic ($E_{k+M}=E_k$ for all $k\in\mathbb N$) and $\sum_{k=1}^M \max E_k\le 0$. Then we can choose real numbers $c_k\ge \max E_k$ ($k=1,\dots, M$) with $\sum_{k=1}^M c_k=0$ and, for $n=\ell M+m$, $0\le m<M$, represent any polynomial $P_n\in\mathcal P$ as
$$
P_n(x)=1+\frac{P(x)(1-x^{\ell M})}{1-x^M}+Q(x)+R(x)\,,
$$
where $P(x)=\sum_{k=1}^M c_k x^k$, $Q(x)=\sum_{k=\ell M+1}^n\eps_k x^k$ and $R(x)$ has non-positive coefficients. Note now that $P(1)=0$, so there is a polynomial $\wt P$ of degree at most $M-1$ such that 
$P(x)=(1-x)\wt P(x)$. Denoting $S(x)=\sum_{k=0}^{M-1}x^k$, we obtain
$$
S(x)P_n(x)=S(x)+\wt P(x)-x^{\ell M}\wt P(x)+S(x)Q(x)+S(x)R(x)\,.
$$
Hence, a coefficient in the expansion of $S(x)P_n(x)$ can be positive only if the corresponding coefficient in the expansion of $S(x)+\wt P(x)-x^{\ell M}\wt P(x)+S(x)Q(x)$ is positive as well.
However, the latter polynomial can have non-zero coefficients only at the powers $k\in[0,M-1]\cup[\ell M+1,(\ell+2)M-2]$. Thus, the coefficient sequence of $S(x)P_n(x)$ can have at most $3M$ sign changes, and, therefore, by the Descartes rule of signs, we have $r\le 3M$ as well. Essentially the same argument with the same conclusion applies to the case $\sum_{k=1}^M\min E_k\ge 0$.

Thus, for $M$-periodic sequences $\mathcal E$, the necessary condition for the possibility to have arbitrarily many roots in $[0,1]$ for some polynomial $P_n\in\mathcal P(\mathcal E)$ is
\begin{equation}
\sum_{k=1}^M \min E_k<0 \quad\text{ and }\quad \sum_{k=1}^M \max E_k>0\,.
\label{pernec}
\end{equation} 
Since we want to deal with not necessarily periodic sequences $\mathcal E$, we generalize \eqref{pernec} as follows.
\begin{definition}
We call a sequence $\mathcal E$ of uniformly bounded closed sets $E_k\subset\R$ {\em balanced} if there exist $M\in\mathbb N$, $a>0$ such that for every $n\ge 0$, one has 
$$
\sum_{k=1}^M \min E_{n+k}\le -a \quad\text{ and }\quad \sum_{k=1}^M \max E_{n+k}\ge a\,.
$$ 
\end{definition}
Note that for periodic sequences the condition of being balanced is equivalent to \eqref{pernec} and that the sequences $\mathcal E$ corresponding to the Littlewood and Newman polynomial families are balanced with parameters $M=2$ and $a=1$, say.

Our main result is the following 
\begin{theorem}
If a sequence $\mathcal E$ of uniformly bounded closed sets $E_k\subset[-A,A]$ is balanced with parameters $M,a$, then for every $r\in\mathbb N$, there exists a polynomial $P_n\in\mathcal P(\mathcal E)$ of degree $n\le C(A,a,M)r^2$ that has at least $r$ distinct roots in $[0,1]$. Moreover, this polynomial can be obtained by an explicit algorithm with running time 
polynomial\,\footnote{The running time bound proved in this paper is $\wt O(r^5)$, where the tilde over $O()$ means that we ignore factors logarithmic in $r$.} in $r$.
\end{theorem}

The rest of the paper is organized as follows. In Section \ref{jensen}, we present the classical Jensen bound on the possible smallness of polynomials from $\mathcal P$ on the interval $I(\alpha)=[1-2\alpha,1-\alpha]$ with $\alpha\in(0,\frac 13)$. In Section \ref{series}, we show how it can be used to force many roots on that interval and reduce the problem to building a power series with restricted coefficients converging to $0$ at finitely many given points. In Section \ref{newman}, we present the Newman Decomposition Lemma (compare with the argument in \cite{BE}, part 3, pages 103-105), which serves as the main tool for all subsequent constructions. Section \ref{trap} reduces the construction of the power series from Section \ref{series} to the investigation of a certain one-dimensional controlled dynamical system. Section \ref{dynamics} is devoted to the analysis of this system and the appropriate control choice. It completes the formal existence proof. Section \ref{algorithm} discusses the corresponding algorithm, its running time, and some details of its implementation. The Appendix contains the proof of the upper bound \eqref{r2len}. 

\section{The Jensen estimate} 
\label{jensen}
\begin{lemma}
Let $\alpha\in(0,\frac 13)$ and let $\mathcal E$ be a uniformly bounded sequence of closed sets $E_k\subset[-A,A]$. Then for every $P_n\in\mathcal P(\mathcal E)$, we have
$$
\int_{I(\alpha)}\log_-|P_n(x)|\,dx\le C(A)\,,
$$
where $I(\alpha)=[1-2\alpha,1-\alpha]$ and $\log_-z=\max(0,-\log z)$.
\label{jenlem}
\end{lemma}
\bp
Consider the polynomal $P_n(z)$ in the domain $\Omega=(\frac 13+\frac 23\mathbb D)\setminus [1-2\alpha,1]$ where $\mathbb D=\{z:|z|<1\}$ is the unit disk. Denote by $\omega$ the harmonic measure on $\partial\Omega$ associated with $0$. We obtain
\begin{equation}
0=\log|P_n(0)|\le\int_{\partial\Omega}\log|P_n|\,d\omega
=\int_{\partial\Omega}\log_+|P_n|\,d\omega-\int_{\partial\Omega}\log_-|P_n|\,d\omega\,.
\label{logs}
\end{equation}
Everywhere in the unit disk, we have the estimate
$$
|P_n(z)|\le 1+A\sum_{k=1}^\infty|z|^k\le\max(1,A)\frac{1}{1-|z|}\,.
$$
Also, on the boundary $\partial\Omega$, we have $\frac{1}{1-|z|}\le \frac 4{|1-z|^2}$. Indeed, on $[1-2\alpha,1]$, we can write 
$$
\frac 1{1-|z|}=\frac 1{|1-z|}\le \frac 1{|1-z|^2}\,,
$$
while on the circle $\frac 13+\frac 23\mathbb T$, we have $|z-\frac 13|^2=\frac 49$, i.e.,
$3|z|^2-2\Re z-1=0$, or
\begin{equation}
|1-z|^2=2(1-|z|^2)=2(1-|z|)(1+|z|)\le 4(1-|z|)\,,
\label{gamma}
\end{equation}
which is equivalent to the claimed inequality.

Thus
\begin{multline*}
\int_{\partial\Omega}\log_+|P_n|\,d\omega\le
\int_{\partial\Omega}\log_+\frac{4\max(1,A)}{|1-z|^2}\,d\omega(z)
\\
=
\int_{\partial\Omega}\log\frac{4\max(1,A)}{|1-z|^2}\,d\omega(z)\le
\log[4\max(1,A)]
\end{multline*}
because $z\mapsto \log\frac{4\max(1,A)}{|1-z|^2}$ is a non-negative harmonic function in $\Omega$, so the integral of its boundary values with respect to $\omega$ does not exceed its value at $0$. Hence, by \eqref{logs},
$$
\int_{\partial\Omega}\log_-|P_n|\,d\omega\le \log[4\max(1,A)]
$$
as well.

It remains to show that $d\omega(z)\ge c|dz|$ on $I(\alpha)$ with some absolute $c>0$. To this end, we first apply the conformal mapping $\zeta=\zeta(z)=-i\frac{z-1}{z+\frac 13}$, which maps $\Omega$ to the upper half-plane with a vertical slit from $0$ to $ih$ with $h=\frac{2\alpha}{\frac 43-2\alpha}$. Note that 
$$
\zeta(0)=3i,\quad
\zeta(1-2\alpha)=ih,\quad \zeta(1-\alpha)=ih'=i\frac{\alpha}{\frac 43-\alpha},\quad \zeta(1)=0\,.
$$
Also,
$$
|d\zeta|=\frac{4}{3|z+\frac 13|^2}|dz|\ge \frac 34|dz|
$$ on $I(\alpha)$. Now apply the second conformal mapping
$\xi=\xi(\zeta)=\sqrt{\zeta^2+h^2}$, which (with an appropriate choice of the branch of the square root) maps the upper half-plane with the slit $[0,ih]$ to the upper half-plane. Then 
$$
\xi(3i)=i\sqrt{9-h^2},\quad \xi([0,ih])=[-h,h]\,,
$$ 
and 
$$|d\xi|=\frac{|\zeta|}{\sqrt{|\zeta^2+h^2|}}|d\zeta|\ge \frac{h'}{h}|d\zeta|$$ 
on $[ih',ih]$. Finally, the harmonic measure on the line $\R$ with respect to the point $i\sqrt{9-h^2}$ is just
$$
\frac 1\pi\frac{\sqrt{9-h^2}}{\xi^2+(9-h^2)}|d\xi|\ge \frac 1\pi\frac{\sqrt{9-h^2}}{9}|d\xi|
$$
on $[-h,h]$.

Bringing all these estimates together and observing that each point $\zeta$ on the slit splits into two points $\xi$ on $\R$, we obtain
$$
d\omega(z)\ge \frac 1\pi\frac{\sqrt{9-h^2}}9\frac{2h'}{h}\frac 34|dz|
$$
for $z\in I(\alpha)$.
It remains to note that $h=\frac{2\alpha}{\frac 43-2\alpha}\le 1$ and $\frac{2h'}{h}=\frac{\frac 43-2\alpha}{\frac 43-\alpha}\ge \frac 23$ for $\alpha\in(0,\frac 13)$, so we can take $c=\frac 1\pi\frac{\sqrt 8}9\frac 23\frac 34=\frac{\sqrt 2}{9\pi}$, say, and finally get
$$
\int_{I(\alpha)}\log_-|P_n(x)|\,dx\le \frac{9\pi}{\sqrt 2}\log[4\max(1,A)]\,.
$$
Of course, we by no means pretend that this bound is sharp.
\ep

\section{Forcing roots}
\label{series}
To ensure that $P_n$ has many roots in $[0,1]$, we will use the following elementary
\bl
Let $I$ be an interval and let $f:I\to\R$ be a continuous function on $I$. Split $I$ into 
$s-1$ equal subintervals $I_j$. If $\frac 1{|I|}\int_I\log_-|f|\le\beta$ and 
$\frac 1{|I_j|}\left|\int_{I_j}f\right|< e^{-2\beta}$ for each $j$, then $f$ has at least $\frac{s-1}2$ roots on $I$. 
\label{manyroots}
\el
\bp
Since $\frac 1{|I|}\int_I\log_-|f|\le\beta$, we can conclude that for at least $\frac{s-1}2$ intervals $I_j$, one has $\frac 1{|I_j|}\int_{I_j}\log_-|f|\le 2\beta$. But for each such $I_j$, we have
\begin{multline*}
\frac 1{|I_j|}\int_{I_j}|f|\ge \exp\left[\frac 1{|I_j|}\int_{I_j}\log|f|\right]
\\
\ge \exp\left[-\frac 1{|I_j|}\int_{I_j}\log_-|f|\right]\ge e^{-2\beta}
>\frac 1{|I_j|}\left|\int_{I_j}f\right|\,,
\end{multline*}
so $f$ has to change sign on $I_j$.
\ep
We shall apply Lemma \ref{manyroots} to the function $f(x)=x^{L-1}P_n(x)$ with some suitably chosen $L>0$, the value $s=2r$, and the interval $I(\alpha)=[1-2\alpha,1-\alpha]$. In this case we shall have by the Jensen estimate (Lemma \ref{jenlem})
$$
\frac 1{|I(\alpha)|}\int_{I(\alpha)}\log_-|f|\le \frac{C(A)}{\alpha}+L\log\frac 1{1-2\alpha}=:\beta\,.
$$
Let now $x_1<x_2<\dots<x_s$ be the endpoints of the $s-1$ intervals of equal length tiling $I(\alpha)$.
Then
$$
\int_{[x_j,x_{j+1}]}f= 
\int_{[x_j,x_{j+1}]}x^{L-1}\left[1+\sum_{k=1}^n\eps_k x^k\right]\,dx
=\left.x^LQ_n(x)\right\vert_{x_j}^{x_{j+1}}\,,
$$
where
$$
Q_n(x)=\frac 1L+\sum_{k=1}^n\frac{\eps_k}{L+k}x^k\,.
$$
Thus, to force $r$ roots of $P_n$ on $I(\alpha)$, it will suffice to ensure that 
$|Q_n(x_j)|< \frac{\alpha}{2(s-1)}e^{-2\beta}$ for all $j=1,\dots,s=2r$.

To this end, we shall construct an infinite series 
$$
Q(x)=\frac 1L+\sum_{k=1}^\infty\frac{\eps_k}{L+k}x^k\,,\quad \eps_k\in E_k\,,
$$
such that $Q(x_j)=0$ for all $j=1,\dots,s$. Truncating it at $k=n$, we will then get
$$
|Q_n(x_j)|\le\sum_{k=n+1}^\infty\frac{|\eps_k|}{L+k}x_j^k\le 
\sum_{k=n+1}^\infty\frac{A}{L+k}(1-\alpha)^k\le \frac A{\alpha(L+n)}e^{-n\alpha}\,.
$$
Bringing all these observations and estimates together, we see that $P_n$ will have at least $r$ roots on $I(\alpha)\subset[0,1]$ if $Q(x_j)=0$ for $j=1,\dots,s$ and $n$ satisfies
\begin{equation}
n\alpha+\log\frac{\alpha(L+n)}A+\log\frac{\alpha}{2(s-1)}>
2\left[\frac{C(A)}{\alpha}+L\log\frac 1{1-2\alpha}\right]\,.
\label{suff}
\end{equation}
Below we will 
prove the following
\begin{proposition}
For every $M\in \mathbb N, A,a>0$, there exists $\eta\in(0,\frac 13)$ and 
$L_0>0$ such that if 
\begin{itemize}
\item a sequence of closed sets $E_k\subset[-A,A]$ is balanced with parameters $M,a$, 
\item $s\in\mathbb N$ and $x_1,\dots,x_s\in(0,1)$ satisfy $\prod_{j=1}^s x_j\ge 1-\eta$, 
\item $L\ge \max(L_0,\frac s\eta)$, 
\end{itemize}
then there exists a power series
$$
Q(x)=\frac 1L+\sum_{k=1}^\infty\frac{\eps_k}{L+k}x^k\,,\quad \eps_k\in E_k\,,
$$
satisfying $Q(x_j)=0$ for all $j=1,\dots,s$.
\end{proposition}
We then can take $s=2r$ and $\alpha=\frac\eta{2s}$. In this case for any $x_j\in I(\alpha)$, we shall have $\prod_{j=1}^s x_j\ge 1-\eta$, so choosing $L=\max(L_0,\frac s\eta)$, we will be able to apply the proposition to establish the existence of the power series $Q(x)$ we need. On the other hand, it is easy to see that, for this choice of parameters, the inequality \eqref{suff} will, indeed, be satisfied for $n=Cr^2$ with sufficiently large $C>0$ depending on $A,a,M$ only. Hence the proof of our theorem will be complete once the proposition is established. 

\section{Newman's Decomposition Lemma}
\label{newman}
\bl
For every $\delta>0$, there exists $\eta=\eta(\delta)\in(0,\frac 13)$ such that for every $s\in\mathbb N$ and every $x_1,\dots,x_s\in(0,1)$ with 
$\prod_{j=1}^s x_j\ge 1-\eta$, there exist $\nu_k\in\R$, $k=0,1,\dots$, satisfying $\sum_{k=0}^\infty|\nu_k|<\delta$ and such that
\begin{equation}
x_j^{-1}-1=\sum_{k=0}^\infty \nu_k \mu ^k x_j^k
\label{newdec}
\end{equation}
for all $j=1,\dots,s$, where $\mu =1-\frac\eta s$\,.
\label{newlem}
\el
\bp
Fix large $\ell=\ell(\delta)\in\mathbb N$ and define
$$
B(z)=\prod_{j=1}^s\frac{1-\mu x_jz}{\mu x_j(\mu x_j-z)}\,,\qquad G_\ell(z)=1-\frac{\ell+1}{\ell}z^{-1}+\frac 1\ell z^{-\ell-1}\,.
$$
For $k=0,1,\dots$, put
$$
c_k=\oint_{\mathbb T}B(z)G_\ell(z)z^k\frac {dz}{2\pi i}\,,
$$
where, as usual, $\mathbb T=\{z:|z|=1\}$ is the unit circle traversed counterclockwise.
Note that
$$
\sum_{k=0}^\infty c_k\mu^kx_j^k=\oint_{\mathbb T}B(z)G_\ell(z)\frac 1{1-\mu x_jz}\frac {dz}{2\pi i}
$$
and the integrand on the right hand side has no poles outside the unit disk because $1-\mu x_jz$ in the denominator cancels with the same factor in the numerator of $B(z)$. Hence, the right hand side evaluates to the residue of the integrand at $\infty$, which is $-\frac 1{\mu x_j}$.

Next we shall find a summable majorant $C_k$ for $c_k$. To this end, we will shift the contour to $\Gamma=\frac 13+\frac 23\mathbb T$. Note that the factor $\frac{1-\mu x_jz}{\mu x_j(\mu x_j-z)}$ maps $\mathbb T$ to $\frac 1{\mu x_j}\mathbb T$ and the circle $\frac 12+\frac 12\mathbb T$ with diameter $[0,1]$ to the circle with diameter $[-\frac 1{\mu x_j},\frac 1{(\mu x_j)^2}]$. The image of $\Gamma$ is squeezed in between, so
$$
|B(z)|\le\prod_{j=1}^s \frac 1{(\mu x_j)^2}\le \frac 1{(1-\eta)^4}\le 6
$$
for every $z\in\Gamma$ as long as $\eta\in(0,\frac 13)$. Now put
$$
C_k=6\int_{\Gamma}|G_\ell(z)||z|^k\frac {|dz|}{2\pi}
$$
and notice that $C_k$ depends neither on $\eta$, nor on the choice of $x_j$.

Observe also that $G_\ell$ has a root of multiplicity $2$ at $1$, so $|G_\ell(z)|\le \gamma_\ell|z-1|^2$ on $\Gamma$ with some $\gamma_\ell\in(0,+\infty)$. Thus
$$
\sum_{k=0}^\infty C_k\le 3\gamma_\ell\int_{\Gamma}\frac{|z-1|^2}{1-|z|}\frac{|dz|}{\pi}
\overset{\eqref{gamma}}\le 12\gamma_\ell \int_{\Gamma}\frac{|dz|}{\pi}=16\gamma_\ell\,.
$$
On the other hand, since $B(z)$ is analytic outside the unit disk, $B(\infty)=1$, 
and $|B(z)|=\prod_{j=1}^s\frac 1{\mu x_j}\le\frac 1{(1-\eta)^2}$ when $z\in\mathbb T$, we have
$$
\int_{\mathbb T}|B(z)-1|^2\frac{|dz|}{2\pi}=\int_{\mathbb T}|B(z)|^2\frac{|dz|}{2\pi}-1
\le\frac{1}{(1-\eta)^4}-1\to 0
$$
as $\eta\to 0$, i.e., $B$ converges to $1$ in $L^2(\mathbb T)$ as $\eta\to 0$. Hence, for every fixed 
$k=0,1,\dots$, we have
$$
c_k\to \oint_{\mathbb T}G_\ell(z)z^k\frac{dz}{2\pi i}=
\begin{cases}
-\frac{\ell+1}{\ell}, &k=0;
\\
\frac 1\ell, &k=\ell;
\\
0 &\text{otherwise}.
\end{cases}
$$
Setting $\nu_0=-\mu c_0-1$, $\nu_k=-\mu c_k$ for $k=1,2,\dots$, we conclude that \eqref{newdec} holds and also
$\sum_{k=0}^\infty |\nu_k|\to \frac 2\ell$  as $\eta\to 0$ by the dominated convergence theorem. It remains to to choose $\ell>\frac 2{\delta}$ and to notice that $c_k$ and, thereby, $\nu_k$ are real because $\vphantom{\overset{a}{A}}B(\bar z)G_\ell(\bar z)=\overline{B(z)}\overline{G_\ell(z)}$.
\ep

\section{Trap $T_{\Psi,\Lambda}$}
\label{trap}
For $x_1,\dots,x_s\in (0,1)$, introduce the notation 
$$
w_k=\begin{bmatrix}
x_1^k
\\
\vdots
\\
x_s^k
\end{bmatrix},\quad k\in\mathbb Z\,.
$$
Then we are looking for a sequence $\eps_k\in E_k$ such that the series 
$\frac 1Lw_0+\sum_{k=1}^\infty \frac{\eps_k}{L+k}w_k$ converges to $0$ in $\R^s$.
Note that since all $x_j\in (0,1)$ and $|\eps_k|\le A$ for all $k$, the series always converges, so it will be enough to ensure that some subsequence of partial sums tends to $0$. Let 
$S=\operatorname{diag}[x_1^{-1},\dots,x_s^{-1}]$ be the linear operator satisfying $Sw_k=w_{k-1}$ for all $k\in\mathbb Z$.

Define 
$$
W(n)=(L+n)\,S^n\!\left[\frac 1L w_0+\sum_{k=1}^n\frac{\eps_k}{L+k}w_k\right]\,.
$$
Then $W(0)=w_0$ and $W(n)$ satisfy the controlled recurrence 
\begin{equation}
W(n+1)=\frac{L+n+1}{L+n}S\,W(n)+\eps_{n+1}w_0\,,
\label{recur}
\end{equation}
where from now on we will think of the mapping $\R^s\ni W\mapsto \frac{L+n+1}{L+n}SW\in \R^s$ as an unstable time dependent dynamical system (the role of time is played by $n$ here) and of the sequence $\eps_n\in E_n$ as of control we can choose to try to stabilize it.

Our goal will be to build a bounded trap $T\subset \R^s$ such that $W(0)=w_0\in T$ and if $W(n)\in T$ for some $n\ge 0$, then we can ensure by an appropriate choice of $\eps_{n+1},\dots,\eps_{n+m}$ that $W(n+m)\in T$ for some $m\ge 1$ again. If such a trap is constructed, then we shall have $W(n)\in T$ or, equivalently,
$$
\frac 1Lw_0+\sum_{k=1}^n\frac{\eps_k}{L+k}w_k\in\frac{1}{L+n}S^{-n}T\,,
$$ 
for infinitely many $n$. But $\frac{1}{L+n}S^{-n}T$ shrinks to $0$ as $n\to\infty$, so the desired convergence to $0$ will be established.

To this end, fix $\delta>0$ to be chosen later and take $\eta=\eta(\delta)\in(0,\frac 13)$ given by the Newman Decomposition Lemma (Lemma \ref{newlem}). Then for any choice of $x_1,\dots,x_s\in (0,1)$ satisfying $\prod_{j=1}^{s}x_j\ge 1-\eta$, we have the decomposition
\begin{equation}
w_{-1}=w_0+\sum_{k=0}^\infty \nu_k\mu^kw_k\quad\text{ with }\sum_{k=0}^\infty|\nu_k|<\delta\,,
\label{wm1}
\end{equation}
where, as before, $\mu =1-\frac{\eta}s$.

We shall be interested in the representations
$$
W(n)=\psi_n w_0+\sum_{k=1}^\infty \lambda_{n,k} \mu ^kw_k,\qquad \psi_n,\lambda_{n,k}\in\R\,.
$$
Since $W(0)=w_0$, we can set $\psi_0=1$, $\lambda_{0,k}=0$ for $k=1,2,\dots$.
Using the recurrence \eqref{recur} and the Newman decomposition \eqref{wm1} of $w_{-1}=Sw_0$,
we see that we can put
\begin{equation}
\psi_{n+1}=\mfrac{L+n+1}{L+n}[\psi_n+\nu_0\psi_n+\mu \lambda_{n,1}]+\eps_{n+1}
\label{dynpsi}
\end{equation}
and
\begin{equation}
\lambda_{n+1,k}=\mfrac{L+n+1}{L+n}[\mu \lambda_{n,k+1}+\psi_n\nu_k]\,.
\label{dynlambda}
\end{equation}
Now put $U_k=\sum_{i=k}^\infty|\nu_i|$, $k=0,1,\dots$, and fix a big $\Lambda>0$.
Assume that $L$ is chosen so that $\frac{L+1}L \mu \le 1$ ($L\ge \frac s\eta$ will suffice for this). Assume also 
that $|\lambda_{n,k}|\le \Lambda U_k$ for all $k=1,2,\dots$ and that $|\psi_n|\le \mu \Lambda$. Note that these inequalities hold for $n=0$, provided that $\Lambda\ge \frac 32$, say.
Then
\begin{multline*}
|\lambda_{n+1,k}|\le\mfrac{L+n+1}{L+n}\bigl[\mu |\lambda_{n,k+1}|+|\psi_n||\nu_k|\bigr]
\\
\le\mfrac{L+1}L\bigl[\mu \Lambda U_{k+1}+\mu \Lambda|\nu_k|\bigr]=\mfrac{L+1}{L}\mu \Lambda U_k\le
\Lambda U_k
\end{multline*}
again. 

Thus, we shall not lose control over $|\lambda_{n,k}|$ before we encounter $\psi_n$ with 
$|\psi_n|>\mu \Lambda$. Our task now is to show that if $\Lambda=\Lambda(A,a,M)$ is chosen large enough and $\delta=\delta(A,a,M)$ is chosen small enough (in this order), then we can keep $|\psi_n|$ under $\mu \Lambda$ for all $n$ by choosing our control $\eps_n$ in an appropriate way. 

From the technical standpoint, it becomes just a question about the controlled one-dimensional dynamics of $\psi_n$ given by \eqref{dynpsi}. It will suffice to show that there exists $\Psi\in[1,\mu \Lambda]$ such that if $|\psi_n|\le\Psi$, then also for some $m\in\mathbb N$,
we have $|\psi_{n+m}|\le\Psi$ while $|\psi_{n+1}|,|\psi_{n+2}|,\dots,|\psi_{n+m-1}|\le \mu \Lambda$ (but not necessarily $\le\Psi$). Then the set
$$
T=T_{\Psi,\Lambda}=
\left\{\psi w_0+\sum_{k=1}^\infty \lambda_k\mu^k w_k\,:\,|\psi|\le\Psi, |\lambda_k|\le\Lambda U_k \right\}
$$
will be our desired trap. Note that all entries of all vectors in $T$ do not exceed 
$\Psi+\Lambda\delta\sum_{k=1}^\infty \mu ^k=\Psi+\frac{\Lambda\delta}{1-\mu }$ in absolute value, so $T$ is, indeed, bounded.

\section{One dimensional controlled dynamics}
\label{dynamics}

Let us see first how fast $\psi_n$ can grow in principle regardless of the choice of the controls $\eps_n$. We have
$$
\psi_{n+1}=\psi_n+\Delta_n+\eps_{n+1}\,,
$$
where
\begin{multline*}
|\Delta_n|=\left|\mfrac 1{L+n}\psi_n+\mfrac{L+n+1}{L+n}(\nu_0\psi_n+\mu \lambda_{n,1})\right|
\\
\le \mfrac 1L \mu \Lambda+\mfrac{L+1}{L}(|\nu_0|\mu \Lambda+\mu \Lambda U_1)
\\
=\mfrac 1L \mu \Lambda+\mfrac{L+1}L \mu \Lambda U_0
\le \left(\mfrac 1L+\delta\right)\Lambda
\end{multline*}
as long as $|\psi_n|\le \mu \Lambda$. We also have $|\eps_{n+1}|\le A$.

Thus, if we start with $|\psi_n|\le\Psi<\mu \Lambda$ and run into trouble after $m$ steps (i.e., have
$|\psi_n|,|\psi_{n+1}|,\dots,|\psi_{n+m-1}|\le \mu \Lambda$ but $|\psi_{n+m}|>\mu \Lambda$),
we must have
$$
m\left[\left(\mfrac 1L+\delta\right)\Lambda+A\right]>\mu \Lambda-\Psi\,.
$$
If we now choose $\Psi=\frac\Lambda 3=3MA+1$ and $\frac 1L,\delta\le \frac 1{9M}$, we shall have
$$
\mu \Lambda-\Psi\ge\frac 23\Lambda-\frac 13\Lambda=\Psi
$$
and
$$
\left(\mfrac 1L+\delta\right)\Lambda+A\le \frac 2{9M}3\Psi+\frac 1{3M}\Psi=\frac\Psi M\,,
$$
so $m>M$. Thus, starting with $|\psi_n|\le\Psi$, we shall be safe with our values of $\psi_{n+1},\dots,\psi_{n+M}$ for any choice of the controls. The choice we will make is the one pushing $\psi_{n+1}$ towards $0$ as hard as possible. More precisely, we will define
\begin{equation}
\eps_{n+1}=
\begin{cases}
\min E_{n+1}, &\psi_n+\Delta_n\ge 0;
\\
\max E_{n+1}, &\psi_n+\Delta_n<0.
\end{cases}
\label{control}
\end{equation}
With this control, two things may happen. 

One possibility is that for some $m\in\{1,\dots,M\}$, there is at least one sign change in the sequence $\psi_{n+m-1}, \psi_{n+m-1}+\Delta_{n+m-1},\psi_{n+m}$. In this case, we have
$$
|\psi_{n+m}|\le |\Delta_{n+m-1}|+|\eps_{n+m}|\le \left(\mfrac 1L+\delta\right)\Lambda+A\le\frac\Psi M\le\Psi\,,
$$
so we have returned to our trap $T$ and may start counting all over.

The other possibility is that all the numbers $\psi_{n+m-1}+\Delta_{n+m-1}$, $m=1,\dots,M$ and $\psi_{n+m}$, $m=0,\dots,M$ are of the same sign and we have $\eps_{n+m}=\min E_{n+m}$ if this sign is positive or $\eps_{n+m}=\max E_{n+m}$ if it is negative all the way from $1$ to $M$. This case requires a bit more careful alalysis. Suppose that the sign is positive and $\eps_{n+m}=\min E_{n+m}$ for 
$m=1,\dots,M$. Then, using the condition that our sequence $\mathcal E$ of sets is balanced with parameters $M,a$, we get
$$
\psi_{n+M}=\psi_n+\sum_{m=1}^M\Delta_{n+m-1}+\sum_{m=1}^M\min E_{n+m}
\le\Psi+M\left(\mfrac 1L+\delta\right)\Lambda-a\,.
$$
Thus, to be certain that $\psi_{n+M}\le\Psi$ in this case, we need to impose yet another condition on $L$ and $\delta$, which is
$$
\frac{1}{L},\delta\le \frac a{2M\Lambda}=\frac a{6M(3MA+1)}\,.
$$
The same condition will suffice for the case when the sign is negative and we choose all maxima.

The final choice of the parameters will then be the following:
\begin{align*}
\Lambda&=3(3MA+1),\quad \Psi=3MA+1,\quad 
\delta=\min\left(\frac 1{9M},\frac a{6M(3MA+1)}\right)\,,
\\
\eta&=\eta(\delta), \quad L_0=\max\left(\frac {6M(3MA+1)}a, 9M\right), 
\quad L\ge\max\left(L_0,\frac s\eta\right)\,.
\end{align*}
It is easy to see that $\Lambda,\Psi,\delta,\eta, L_0$ depend on $A,a,M$ only, so the proof of the proposition is complete.

\section{The algorithm}
\label{algorithm}

It should be clear now how to build an algorithm for finding a polynomial $P_n\in \mathcal P(\mathcal E)$ with a given number $r$ of roots in principle. Given $r$, one should choose $s\ge r+1$, $n\asymp r^2$, $L\asymp r$, $\eta\asymp 1$ and set $\alpha=\frac{\eta}{2s}$, $\mu=1-\frac\eta s$. Then one needs to choose the points $x_1,\dots,x_s\in I(\alpha)$ and compute $\nu_0,\nu_1,\dots,\nu_n$ (the rest of the coefficients in the Newman decomposition are never used in the determination of $\eps_1,\dots,\eps_n$). Finally, one can initialize $\psi=1,\lambda_k=0$ and update their values according to \eqref{dynpsi} and \eqref{dynlambda} choosing the coefficient sequence $\eps_n$ as in \eqref{control}. 

The total running time is the time of the pre-computation of $\nu_k$ plus the time of the determination of the coefficients. Since each coefficient determination requires updating an array of length $n$, we get $n$ steps that take $\asymp n$ elementary arithmetic operations each. So it looks like we need about $n^2\asymp r^4$ operations for this part. However, as it is usual with the real number computations, the question of the propagation of the rounding errors arises.

In theory, we are just running the exponentially unstable dynamics
$$
W(0)=w_0,\qquad W(m+1)=\frac{L+m+1}{L+m}SW(m)+\eps_{m+1}w_0
$$ 
and trying to stabilize it by choosing an appropriate control sequence $\eps_m$. A small issue, however, is that we determine $\eps_m$ not from the vectors $W(m)$ directly, but rather from their representations
$$
W(m)=\psi_m w_0+\sum_{k=1}^\infty \lambda_{m,k}\mu ^kw_k\,,
$$
so the adequate computation model seems to be given by
\begin{equation}
\begin{aligned}
\wh W(m)&=\wh\psi_m w_0+\sum_{k=1}^\infty \wh\lambda_{m,k}\mu ^kw_k\,,
\\
\wh\psi_{m+1}&=\mfrac{L+m+1}{L+m}(\wh\psi_m+\nu_0\wh\psi_m+\mu \wh\lambda_{m,1})+\wh\eps_{m+1}+O(\tau)\,,
\\
\wh\lambda_{m+1,k}&=\mfrac{L+m+1}{L+m}(\mu \wh\lambda_{m,k+1}+\wh\psi_m\nu_k)+O(\tau)\,,
\end{aligned}
\label{erratic}
\end{equation}
where $\tau$ is the fixed point rounding error and $\wh\eps_m$ are determined as in $\eqref{control}$ but using the (erratic) values $\wh\psi_m$ and $\wh\lambda_{m,k}$ instead of the true $\psi_m$ and $\lambda_{m,k}$, i.e,
\begin{equation*}
\wh\eps_{m+1}=
\begin{cases}
\min E_{m+1}, &\mfrac{L+m+1}{L+m}(\wh\psi_m+\nu_0\wh\psi_m+\mu \wh\lambda_{m,1})\ge 0;
\\[7pt]
\max E_{m+1}, &\mfrac{L+m+1}{L+m}(\wh\psi_m+\nu_0\wh\psi_m+\mu \wh\lambda_{m,1})<0.
\end{cases}
\end{equation*}
It is not hard to see that if we allow ourselves some extra leeway in the inequalities for $\Lambda$, $\Psi$, $L$ and $\delta$, then we will not need very high precision to stabilize the erratic dynamics. Indeed, as far as $\wh\lambda_{m,k}$ are concerned, we just notice that if 
$|\wh\lambda_{m,k}|\le\delta+\Lambda U_k$ and $\wh\psi_m\le \mu \Lambda$, we can write
\begin{multline*}
|\wh\lambda_{m+1,k}|\le \mfrac{L+1}{L}[\mu (\delta+\Lambda U_{k+1})+\mu \Lambda\nu_k]+O(\tau)
\\
\le [\mfrac{L+1}{L}\mu \delta+O(\tau)]+\Lambda U_k\le\delta+\Lambda U_k\,,
\end{multline*}
as long as $\mfrac{L+1}{L}\mu \delta<\delta$ and $O(\tau)$ is too small to span the difference. Recalling that 
$\mu =1-\frac{\eta}{s}$ and taking $L\ge\frac{2s}\eta$, we see that the precision $\tau\asymp r^{-1}$ is already enough to keep the values $\wh\lambda_{m,k}$ bounded by $2\delta$ as long as $|\wh\psi_m|$ remain bounded by $\mu \Lambda$. Stabilizing the one-dimensional dynamics of $\wh\psi_m$ then imposes an even weaker restriction on $\tau$. We just need to add an extra $O(\tau)$ term to the estimate for $|\Delta_m|$ in all previous calculations, where we can afford even a constant leeway.

However, the real issue is not the stabilization of the erratic dynamics per se, but ensuring that the true dynamics with the controls $\wh\eps_m$ based on the erratic computations results in not too large vectors $W(m)$. More precisely, we need to show that for the sequence of vectors defined by
$$
W(0)=w_0,\qquad W(m+1)=\frac{L+m+1}{L+m}SW(m)+\wh\eps_{m+1}w_0\,,
$$
we still have a decent bound on the size of $W(n)$. To this end, we will just compare $W(m)$ to $\wh W(m)$, for which we know a good bound from the bounds on $\wh\psi_m$ and $\wh\lambda_{m,k}$. Taking \eqref{erratic} into account,
we see that $\wh W(m)$ satisfy
$$
\wh W(0)=w_0, \qquad \wh W(m+1)=\frac{L+m+1}{L+m}S\wh W(m)+\wh\eps_{m+1}w_0+ \frac 1{1-\mu }O(\tau)\,,
$$  
where the last term is obtained by summing up the rounding errors in $\wh\psi_{m+1}$ and $\wh\lambda_{m+1,k}$, i.e., evaluating the sum 
$\sum_{k=0}^\infty \mu ^k O(\tau)$.

Thus,
$$
W(m+1)-\wh W(m+1)=\frac{L+m+1}{L+m}S(W(m)-\wh W(m))+ \frac 1{1-\mu }O(\tau)\,,
$$
whence, by induction on $m$,
$$
\bigl|W(m)-\wh W(m)\bigr|\le \frac 1{1-\mu }O(\tau)\sum_{k=0}^{m-1}\|S\|^k\le \frac 1{1-\mu }O(\tau)\frac{(1-2\alpha)^{-m+1}}{2\alpha} 
$$
because $\|S\|\le\frac{1}{1-2\alpha}$ as all $x_j\ge 1-2\alpha$. Recalling that $\mu =1-\frac \eta s$ and $\alpha=\frac\eta {2s}$, we see that
$$
\bigl|W(n)-\wh W(n)\bigr|\le O(s^2)e^{O(n/s)}O(\tau)\,.
$$
Since $s\asymp r$ and $n\asymp r^2$, we conclude that to keep the difference $W(n)-\wh W(n)$ bounded in this computational model, we need to choose $\tau=e^{-Cr}$ with sufficiently large $C$.

Now we can address the question on the required precision of the computation of $\mu $ and $\nu_k$. It should be clear at this point that to have our model justified, their values must be computed with the same precision $\tau$. That is not a big deal for $\mu $, which is given by a simple arithmetic formula in terms of our parameters but the computation of $\nu_k$ that are obtained by the contour integration requires a separate discussion.
Recall that $\nu_k$ are obtained by elementary expressions from
$$
c_k=\oint_{\mathbb T}B(z)G_\ell(z)z^k\frac{dz}{2\pi i}\,,
$$
so it will suffice to compute $c_k$ with precision $\tau$ for $k=0,\dots,n$. To this end, we suggest just to discretize the integral to the sum
$$
\frac 1N\sum_{j=0}^{N-1}B\left(e^{2\pi i j/N}\right)G_\ell\left(e^{2\pi i j/N}\right)e^{2\pi i j(k+1)/N}
$$ 
for sufficiently large $N>n+1$ that is a power of $2$ and to use the fast Fourier transform. Note that, for $0\le k<N-1$, this sum, even if computed exactly, is {\em not} the true value of $c_k$ but $c_k+c_{k+N}+c_{k+2N}+\dots$. Fortunately, since all poles of $B(z)$ and $G_\ell(z)$ lie reasonably deep inside the unit disk, $|c_k|$ decays fairly fast as $k\to\infty$. To get a simple but sufficient for our purposes bound, we will just shift the contour to $u\mathbb T$ with $u=1-\frac\eta{2s}$ and observe that for each Blaschke factor $\frac{1-\mu x_jz}{\mu x_j(\mu x_j-z)}$, we have
$$
|1-\mu x_juz| \le|1-u|+|1-\mu x_jz|\le\mfrac 32|1-\mu x_jz|\,,
$$
while
$$
|\mu x_j-uz|\ge |\mu x_j-z|-|1-u|\ge\mfrac 12|\mu x_j-z|\,,
$$
so $|B(uz)|\le 3^s|B(z)|\le \frac {3^s}{(1-\eta)^2}$ when $z\in\mathbb T$.
As to $G_\ell$, we just estimate $\max_{z\in\mathbb T}|G_\ell(uz)|\le u^{-\ell-1}\max_{z\in\mathbb T}|G_\ell(z)|$. At last, $u^k=(1-\frac\eta{2s})^k\le e^{-\frac{\eta k}{2s}}$. Thus, to make the sum $c_{k+N}+c_{k+2N}+\dots$ negligible (i.e., less than $\tau=e^{-Cr}$), we can take any $N$ for which $3^se^{-\frac{\eta N}{2s}}$ is substantially smaller than $e^{-Cr}$, which forces us to choose $N\asymp r^2$ too. 

Thus, the cost of the pre-computation is about $sN+N\log N\asymp r^3$ elementary arithmetic operations, which is $r$ times less than the cost of computing all the coefficients $\wh\eps_m$, $m=1,\dots,n$. The total running time is then about $r^4$ times the time needed for an elementary arithmetic operation on $Cr$-digit numbers, which is $\wt O(r)$. That gives $\wt O(r^5)$ claimed in the statement of the theorem.

Finally, a few words about the practical implementation, if somebody feels a desire to try it. While the orders of magnitude in the above discussion are all correct, the numerical constants given by our rigorous proofs are certainly suboptimal, so the best way to choose an appropriate value for $\eta$ and $\ell$ is by trial and error. This won't waste too much time because if the blow-up in our dynamical system occurs at all, it usually happens rather fast and can be seen after about $r$ iterations already. Also, while the theory guarantees $r$ roots for $s=2r$, in practice $s=r+1$ may already be enough. 

 
\section*{Appendix}

In this section we shall prove the classical bound $r\le C(A)\sqrt n$ for the number $r$ of the roots of a polynomial $P_n(x)=1+\sum_{k=1}^n\eps_k x^k$ with $|\eps_k|\le A$ on the interval $[0,1]$. This bound holds for the number of roots counted with multiplicity, so no assumption that the roots are distinct will be required in the proof.
Since $P_n(0)=1$, it is enough to get an estimate for the number of roots of $P_n$ in $(0,1]$.
We shall follow the exposition in \cite{B} (Theorem 5 on page 55) with some minor modifications.

Suppose that one can construct a polynomial $q$ of degree $m$ with real roots such that
$$
q(0)=1,\qquad \sum_{k=1}^n|q(k)|<\frac 1A\,.
$$
Then the polynomial $\wt P_n(x)=1+\sum_{k=1}^n q(k)\eps_k x^k$ will have no roots in $(0,1]$ because the constant term $1$ dominates the sum of the absolute values of all other terms. However, writing $q(x)=\gamma\prod_{j=1}^m (x-\rho_j)$ with $\gamma,\rho_j\in\R$, we can express $\wt P_n$ as
$$
\wt P_n=\gamma\Bigl[\prod_{j=1}^m D_{\rho_j}\Bigr] P_n\,,
$$
where $D_\rho f=(x\frac d{dx}-\rho)f$.

Since $D_\rho f = x^{\rho+1}\frac d{dx}(x^{-\rho}f)$, each application of $D_{\rho}$ diminishes the number of roots of a function $f$ on the interval $(0,1]$ by at most $1$. Indeed, Rolle's theorem guarantees the existence of a root between two distinct roots, and the multiplicity of a repeated zero drops by $1$. Since after $m$ such applications the polynomial $P_n$ loses all its roots on that interval, we conclude that its initial number of roots $r$ satisfied $r\le m$.

Now, to construct the polynomial $q$ of low degree with the desired property, for $\ell\in\mathbb N$, consider the normalized Dirichlet kernel
$$
\frac 1{2\ell+1}\left[1+2\sum_{k=1}^\ell \cos(ky)\right]
=\frac 1{2\ell+1}\frac{\sin(\ell+\frac 12)y}{\sin\frac y2}\,.
$$
It can be written as $q_0(\cos y)$ where $q_0$ is a polynomial of degree $\ell$ having $\ell$ roots on $[-1,1]$. We have $q_0(1)=1$ and $|q_0(t)|\le \frac 1{2\ell+1}\sqrt{\frac 2{1-t}}$ for $t\in[-1,1)$. Now put $q_1(t)=q_0(1-\frac{2t}n)$. Then $q_1$ is also of degree $\ell$, still has all its roots real, $q_1(0)=1$, and $|q_1(k)|\le \frac 1{2\ell+1}\sqrt{\frac nk}$ for $k=1,\dots,n$. Taking $\ell=\lceil \sqrt n\,\rceil$, we conclude that $q_1(k)\le \frac 1{2\sqrt k}$ for $k=1,\dots,n$. But then for every integer power $v\ge 4$, we have
$$
\sum_{k=1}^n |q_1(k)|^v\le 2^{-v}\sum_{k=1}^\infty\frac 1{k^2}\le 2^{1-v}\,,
$$
so $q=q_1^v$ with some $v=v(A)$ will satisfy the desired property and have the degree 
$m=v(A)\lceil \sqrt n\,\rceil$, yielding the claimed bound for $r$.


\begin{thebibliography}{3}

\bibitem[B]{B}
{\sc P. Borwein}, {\em Computational Excursions in Analysis and Number Theory},
Springer-Verlag, New York, 2002.
\bibitem[BE]{BE}
{\sc P. Borwein and T. Erd\'elyi}, {\em The $L_p$ version of Newman's inequality for lacunary polynomials}, Proceedings of the AMS,
v. 124, no. 1, 1996, pp. 101-109.
\bibitem[BEK]{BEK}
{\sc P. Borwein, T. Erd\'elyi, and G. K\'os}, {\em Littlewood-type problems on $[0,1]$}, Proceedings of the London Math. Soc. v. 79, issue 1, 1999, pp. 22-46.
\bibitem[E]{E}
{\sc T. Erd\'elyi}, {\em Extensions of the Bloch-P\'olya theorem on the number of real zeros of polynomials}, Journal de th\'eorie des nombres de Bordeaux, v. 20, no. 2, 2008, pp. 281-287.
\bibitem[MOP]{MOP}
https://mathoverflow.net/questions/461631/number-of-real-roots-of-0-1-polynomial.

\end{thebibliography}
\end{document}